\documentclass[a4paper]{article}
\usepackage{amsmath, amssymb, eucal, amscd}

\newtheorem{thm}{Theorem}
\newtheorem{lem}[thm]{Lemma}
\newtheorem{eg}[thm]{Example}
\newtheorem{prop}[thm]{Proposition}
\newtheorem{cor}[thm]{Corollary}
\newtheorem{defn}[thm]{Definition}
\newtheorem{rem}[thm]{Remark}

\newtheorem{rem-eg}[thm]{Remark and Example}

\newcommand{\id}{{\rm id}}
\newcommand{\ind}{{\rm ind}}
\newcommand{\I}{{\rm Int}}
\newcommand{\Ind}{{\rm Ind}}
\newcommand{\Res}{{\rm Res}}

\newcommand{\ti}{\tilde}
\newcommand{\la}{\langle}
\newcommand{\ra}{\rangle}
\newcommand{\smnoind}{\smallskip\noindent}

\newcommand{\pr}{{\rm Prim}}
\newcommand{\pq}{{\rm Qorb}}
\newcommand{\pf}{{\rm Prim}_\beta}

\newcommand{\po}{{\rm Orb}}
\newcommand{\abs}[1]{\left\vert#1\right\vert}

\begin{document}
\title{Invariant ideals of twisted crossed products}
\author{Chi-Wai Leung and Chi-Keung Ng\footnote{This work is partially supported by Hong Kong RGC Direct Grant}}
\maketitle
\begin{abstract}
We will prove a result concerning the inclusion of non-trivial invariant ideals inside non-trivial ideals of a twisted crossed product. 
We will also give results concerning the primeness  and simplicity of crossed products of twisted actions of locally compact groups on $C^*$-algebras. 
\end{abstract}

\bigskip

\bigskip

The motivation of this study is our recent research on $C^*$-unique groups (or \emph{Fourier groups}) in \cite{LN}. 
In fact, one interesting question is when the semi-direct product of a $C^*$-unique group with another group is again $C^*$-unique. 
This turns out to be related to the following question: 
\emph{Given a $C^*$-dynamical system $(A, G, \alpha)$, under what condition will it be true that any non-zero ideal of $A\times_\alpha G$ contains a non-zero $\hat\alpha$-invariant ideal ($\hat\alpha$ being the dual coaction)?} 
The aim of this paper is to study this question. 

\medskip

In fact, in the case of discrete amenable groups acting on compact spaces, Kawamura and Tomiyama gave (in \cite{KT}) a complete solution of the above question. 
In \cite{AS}, Archbold and Spielberg generalised the main result in \cite{KT} to the case of discrete $C^*$-dynamical system. 
In this article, we are going to present a weaker result but in the case of general locally compact groups. 
As a corollary, we obtain some equivalent conditions for the primeness of crossed products (in terms of the actions). 
Moreover, we will also give a brief discussion on the simplicity of crossed products (which is related but does not need the main theorem). 

\medskip

Throughout this article, $A$ is a $C^*$-algebra, $G$ is a locally compact group and $(\beta, u)$ is a twisted action of $G$ on $A$ (in the sense of Busby and Smith; see \cite[2.1]{PR}).

\medskip

\begin{rem}
Suppose that $A$ is separable and $G$ is second countable. 

\smnoind
(a) By the stabilisation trick of Packer and Raeburn (\cite[3.4]{PR}), there is a canonical continuous action of $G$ on $\pr(A) \cong \pr(A\otimes \mathcal{K}(L^2(G)))$ which is the same as the canonical one defined by $\beta$.
This means that the canonical action of $G$ on $\pr(A)$ is continuous. 
For any $I\in \pr(A)$, the \emph{quasi-orbit, $\mathcal{Q}(I)$, of $I$} is defined to be the set $\{ I'\in \pr(A): \overline{G\cdot I} = \overline{G\cdot I'}\}$. 
Moreover, $(\beta,u)$ is said to be \emph{free} if $\pr(A)^t=\emptyset$ for any $t\neq e$. 

\smnoind
(b) We recall the following materials from \cite[\S 2]{PR2}. 
An ideal $I$ of $A$ is said to be \emph{$\beta$-invariant} if $\beta_t(I)\subseteq I$ for any $t\in G$. 
Note that $I$ is $\beta$-invariant if and only if $I\otimes \mathcal{K}(L^2(G))$ is invariant under the (ordinary) action given by \cite[3.4]{PR}. 
If $I$ is $\beta$-invariant, we can define a canonical twisted $C^*$-dynamical system $(I,G,\beta\!\mid_I, u\!\mid_I)$.
Moreover, $I\times_{\beta, u} G$ can be regarded as an ideal of $A\times_{\beta, u} G$. 
\end{rem}

\medskip

Let us start our discussion with two simple lemmas. 
The first one is probably well known (although we cannot find it in the literatures). 
As usual, we will equip the orbit space $\po(\beta) = \pr(A)/G$ and the quasi-orbit space $\pq(\beta) = (\pr(A)/G)\ti{\ }$ with the quotient topologies. 

\medskip

\begin{lem}
\label{inv-open}
(a) Suppose that $q_G^{\ }$ is the canonical quotient map from $\pr(A)$ to $\pq(\beta)$. 
An open (or a closed) subset $E$ of $\pr(A)$ is $G$-invariant if and only if it is \emph{quasi-invariant} in the sense that $E = q_G^{-1}(q_G^{\ }(E))$. 

\smnoind
(b) The canonical quotient map $\psi$ from $\po(\beta)$ to $\pq(\beta)$ gives a bijective correspondence between Borel sets in $\po(\beta)$ and Borel sets in $\pq(\beta)$. 
Moreover, if $\gamma$ is any ordinal, $\psi$ will induce a bijection between the $\mathcal{G}_\gamma$ (respectively, $\mathcal{F}_\gamma$) sets in $\po(\beta)$ and those in $\pq(\beta)$ (where $\mathcal{G}_0$ is the collection of all open sets, $\mathcal{G}_1 = \mathcal{G}_\delta$, $\mathcal{G}_2 = \mathcal{G}_{\delta \sigma},...$ and $\mathcal{F}_0$ is the collection of all closed sets, $\mathcal{F}_1 = \mathcal{F}_\sigma$, $\mathcal{F}_2 = \mathcal{F}_{\sigma \delta},...$).
\end{lem}
\noindent
{\bf Proof:}
(a) It is clear that any quasi-invariant subset is $G$-invariant. 
Moreover, a subset $E$ is quasi-invariant if and only if $\mathcal{Q}(I)$ is a subset of $E$ for any $I\in E$.
Suppose that $E$ is $G$-invariant and $I\in E$. 
If $E$ is closed, then 
$$\mathcal{Q}(I)\subseteq \overline{G\cdot I} \subseteq E$$ 
and $E$ is quasi-invariant. 
On the other hand, suppose that $E$ is open.  
For any $I'\in \mathcal{Q}(I)$, we have $I\in \overline{G\cdot I'}$ and so $E\cap G\cdot I' \neq \emptyset$.
Thus, there exists $r\in G$ such that $I'\in r\cdot E\subseteq E$. 

\smnoind
(b) Let $\mathcal{B} = \{ X\subseteq \po(\beta): \psi(X)$ is a Borel set and $\psi^{-1}(\psi(X)) = X\}$.
By part (a), $\mathcal{G}_0 \cup \mathcal{F}_0 \subseteq \mathcal{B}$. 
For any $X\in \mathcal{B}$, we have 
$$\psi^{-1}(\pq(\beta)\setminus \psi(X)) = \po(\beta)\setminus \psi^{-1}(\psi(X)) = \po(\beta)\setminus X.$$ 
Therefore, $\psi(\po(\beta)\setminus X) = \pq(\beta)\setminus \psi(X)$ (as $\psi$ is surjective) is a Borel set and 
$$\psi^{-1}(\psi(\po(\beta)\setminus X)) = \po(\beta)\setminus X.$$ 
Moreover, for any sequence $X_n\in \mathcal{B}$, it is clear that $\psi(\bigcup_{n=1}^\infty X_n) = \bigcup_{n=1}^\infty \psi(X_n)$ is a Borel set and
$$\psi^{-1}(\psi(\bigcup_{n=1}^\infty X_n)) = \bigcup_{n=1}^\infty \psi^{-1}(\psi(X_n)) = \bigcup_{n=1}^\infty X_n.$$ 
Therefore, $\mathcal{B}$ contains all the Borel sets of $\po(\beta)$ and $X\mapsto \psi(X)$ give an injective correspondence into the Borel sets of $\pq(\beta)$. 
On the other hand, since $\psi^{-1}(Y)$ is a Borel set for any Borel set $Y\subseteq \pq(\beta)$ (as $\psi$ is continuous), this correspondence is bijective. 
The second statement follows from similar arguments as above and transfinite inductions. 

\bigskip

Consequently, if every singleton of the orbit space is a Borel set, then all $G$-orbits are quasi-orbits. 

\medskip

The next lemma is crucial for the proof of the main theorem (Theorem \ref{twisted-case}).
Note that amenability of $G$ is required because we need to use the results in \cite{GL} and \cite{GR}. 
From now on, we denote $\pf(A) = \bigcap_{t\in G\setminus \{e\}} \pr(A)\setminus \pr(A)^t$.

\medskip

\begin{lem}
\label{unique}
Suppose that $G$ is a second countable amenable group and $A$ is a separable $C^*$-algebra with an (ordinary) action $\alpha$ by $G$. 
Let $X$ be a non-empty quasi-invariant subset of ${\rm Prim}_\alpha (A)$. 
For any $I\in X$, there exists exactly one $J\in \pr(A\times_\alpha G)$ such that $J$ lives on $\mathcal{Q}(I)$ (in the sense of \cite[p.221]{Gre}). 
In particular, we can take $X=\I({\rm Prim}_\alpha (A))$ if it is non-empty. 
\end{lem}
\noindent
{\bf Proof:}
As $X$ is quasi-invariant, $\mathcal{Q}(I) \subseteq X$ and so for any $I'\in \mathcal{Q}(I)$, 
$$G\!_{I'} = \{e\}$$ 
where $G\!_{I'}$ is the stabilizer of $I'$ (note that $X \cap \pr(A)^t = \emptyset$ for any $t\neq e$). 
On the other hand, as the quasi-orbit map (see e.g. \cite[p.620]{GL}) from $\pr(A\times_\alpha G)$ to $\pq(\alpha)$ is surjective (see e.g. \cite[4.8]{GL}), there exists $J \in \pr(A\times_\alpha G)$ such that $J$ lives on $\mathcal{Q}(I)$. 
Suppose that $J_1, J_2\in \pr (A\times_\alpha G)$ both live on $\mathcal{Q}(I)$. 
By \cite[3.2]{GR}, we see that there exist $I_1, I_2\in \mathcal{Q}(I)$ and representations $\mu_k = \pi_k\times u_k$ of $A\times_\alpha G\!_{I_k} = A$ such that 
$$\ker (\pi_k) \,= \, I_k \qquad {\rm and} \qquad J_k \, = \, \ker (\ind_{(e)}^G (\mu_k)) \, = \, \ker (\Ind (\mu_k))$$ 
for $k=1,2$ (where $\Ind(\mu_k)$ is the regular representation induced by $\mu_k$ as in \cite[p.602]{GL}). 
As $I_1,I_2\in \mathcal{Q}(I)$, we have, 
$$\Res(\Ind{\ }I_1) = \ker (\mathcal{Q}(I)) = \Res(\Ind{\ }I_2)$$ 
(where $\Res = \Res^G_{(e)}$; see \cite[\S 3]{GL}). 
Therefore, by \cite[3.11(i)\&(v)]{GL}, $J_1 = \Ind{\ }I_1 = \Ind(\Res(\Ind{\ }I_1)) = \Ind(\Res(\Ind{\ }I_2)) = \Ind{\ }I_2 = J_2$ as required. 
Finally, $\I({\rm Prim}_\alpha (A))$ is quasi-invariant by Lemma \ref{inv-open}(a) (note that it is $G$-invariant as ${\rm Prim}_\alpha (A)$ is).

\bigskip

Inspired by the condition in \cite{AS}, we make the following definition. 

\medskip

\begin{defn}
A twisted action $(\beta,u)$ is said to be \emph{almost-free} if $\bigcup_{t\in G\setminus \{e\}} \pr(A)^t$ is nowhere dense in $\pr(A)$. 
\end{defn}

\medskip

If $q_G^{\ }$ is the map as in Lemma \ref{inv-open}, it is natural to ask whether there is any relation between the nowhere density of $q_G^{\ }(\bigcup_{t\in G\setminus \{e\}} \pr(A)^t)$ and that of $\bigcup_{t\in G\setminus \{e\}} \pr(A)^t$. 
In the following, we will show that they are in fact the same and are equivalent to some other interesting conditions. 

\medskip

\begin{prop}
\label{almost-free}
For a separable twisted $C^*$-dynamical system $(A,G,\beta,u)$, the following statements are equivalent. 

\smnoind
(i) $q_G^{\ }(\bigcup_{t\in G\setminus \{e\}} \pr(A)^t)$ is nowhere dense in $\pq(\beta)$.

\smnoind
(ii) $q_G^{\ }(\I(\pf(A))$ is dense in $\pq(\beta)$.

\smnoind
(iii) $(\beta,u)$ is almost-free.

\smnoind
(iv) $\I(U\cap\pf(A)) \neq \emptyset$ for any non-empty $G$-invariant open set $U$ of $\pr(A)$. 

\smnoind
(v) The action of $G$ is free on an open dense $G$-invariant subset $U\subseteq \pr(A)$. 

\smnoind
(vi) There is a $\beta$-invariant ideal $I$ of $A$ such that $\ker \{J\in \pr(A): I\nsubseteq J\} = (0)$ and the restriction of $(\beta, u)$ on $I$ is free. 
\end{prop}
\noindent
{\bf Proof:}
Let $F$ be a $G$-invariant subset of $\pr(A)$ and $F^c = \pr(A)\setminus F$.  
It is clear that $q_G^{-1}(\pq(\beta)\setminus q_G^{\ }(F^c)) \subseteq F$ and hence $\I(\pq(\beta)\setminus q_G^{\ }(F^c)) \subseteq q_G^{\ }(\I(F))$ (because $q_G^{\ }$ is continuous and surjective). 
Moreover, we actually have
\begin{equation}
\label{int}
q_G^{\ }(\I(F)) = \I(\pq(\beta)\setminus q_G^{\ }(F^c)).
\end{equation}
In fact, suppose there exists $\mathcal{Q}\in q_G^{\ }(\I(F))\cap q_G^{\ }(F^c)$. 
Take $I\in \I(F)$ and $I'\in F^c$ such that $\mathcal{Q}(I) = \mathcal{Q} = \mathcal{Q}(I')$. 
As $F$ is $G$-invariant, $\I(F)$ is also $G$-invariant and hence quasi-invariant (by Lemma \ref{inv-open}(a)). 
Thus, $I'\in \mathcal{Q}(I)\subseteq \I(F)$ which contradicts $I'\in F^c$. 
As $q_G^{\ }$ is open, we have $q_G^{\ }(\I(F)) \subseteq \I(\pq(\beta)\setminus q_G^{\ }(F^c))$. 

\smnoind
(i)$\Leftrightarrow$(ii) By putting $F = \pf(A)$ into (\ref{int}), we have
$$q_G^{\ }(\I(\pf(A))) = \I(\pq(\beta)\setminus q_G^{\ }(\bigcup_{t\in G\setminus \{e\}} \pr(A)^t))$$ 
and the equivalence is clear.

\smnoind
(ii)$\Rightarrow$(iii) For any non-zero ideal $I_1$ of $A$, we have 
$$q_G^{\ }(\I(\pf(A)))\cap q_G^{\ }(\{ I\in \pr(A): I_1\nsubseteq I\}) \neq \emptyset.$$ 
Thus, there exists $I\in \I(\pf(A))$ and $I'\in \pr(A)$ such that $I_1\nsubseteq I'$ and $\mathcal{Q}(I) = \mathcal{Q}(I')$. 
As $\I(\pf(A))$ is quasi-invariant, $I'\in \mathcal{Q}(I) \subseteq \I(\pf(A))$. 
Therefore, 
$$\I(\pf(A)) \cap \{ I\in \pr(A): I_1\nsubseteq I\}\neq \emptyset$$ 
for any ideal $I_1$. 
This means that $\I(\pf(A))$ is dense. 

\smnoind
(iii)$\Rightarrow$(iv) If $\pf(A)$ ($= \pr(A)\setminus \bigcup_{t\in G\setminus \{e\}} \pr(A)^t$) is everywhere dense, then condition (iv) clearly holds. 

\smnoind
(iv)$\Rightarrow$(ii)
Let $V$ be a non-empty open subset of $\pq(\beta)$ and $U = q_G^{-1}(V)$. 
It is clear that $q_G^{\ }(I) \in V\cap q_G^{\ }(\I(\pf(A)))$ for any $I\in \I(U\cap \pf(A))$ and this establishes condition (ii). 

\smnoind
(iii)$\Leftrightarrow$(v) It is clear that the action of $G$ on an open subset $U\subseteq \pr(A)$ is free if and only if $U \subseteq \I(\pf(A))$.

\smnoind
(v)$\Rightarrow$(vi) Let 
$I = \ker (\pr(A)\setminus U)$. 
As $U$ is a $G$-invariant dense subset, $I$ is $\beta$-invariant and $\ker \{J\in \pr(A): I\nsubseteq J\} = (0)$. 
It is well known that the canonical ($G$-covariant) map from $U$ to $\pr(I)$ is a homeomorphism (see e.g. \cite[5.5.5]{Mur}). 
Therefore, the restriction of $(\beta, u)$ on $I$ is free. 

\smnoind
(vi)$\Rightarrow$(iii) Let $V$ be the dense open set $\{J\in \pr(A): I\nsubseteq J\}$. As the action of $G$ on $V \cong \pr(I)$ is free, $V\subseteq \pf(A)$. 

\medskip

\begin{rem}
(a) In the case of discrete group actions, almost-freeness is stronger than the \emph{topological freeness} in \cite{AS}: $\bigcap_{i=1}^n \hat A\setminus \hat A^{t_i}$ is dense in $\hat A$ for all $t_1,...,t_n\in G\setminus \{e\}$ (see Example \ref{non-eg} below).

\smnoind
(b) When $G=\mathbb{Z}$, almost-freeness (in the form of condition (v) above) is stronger than the requirement that the Connes' spectrum is the whole torus (see \cite[10.4(v)]{OP3}). 
Again, Example \ref{non-eg}(b) gives an action whose Connes's spectrum is the whole torus but is not almost-free. 
\end{rem}

\medskip

\begin{eg}
\label{eg}
(a) Let $X = \mathbb{R}$ with addition and the usual topology and $G = \mathbb{R}^+$ with multiplication and the usual topology. 
Consider the action of $G$ on $X$ by $\alpha_t(n) = tn$ ($t\in \mathbb{R}^+$; $n\in \mathbb{R}$). 
Then $X^t = \{0\}$ for all $t\neq 1$ and so the corresponding action of $G$ on $C_0(X)$ is almost-free. 

\smnoind
(b) Let $X = \mathbb{R}^{n+1}$ and $G = \mathbb{R}^n$ both with additions and the usual topologies. 
Define an action $\alpha$ of $G$ on $X$ by $\alpha_t(u,a) = (u-at,a)$ ($t,u\in \mathbb{R}^n$; $a\in \mathbb{R}$). 
Then for any $t\neq 0$, $X^t = \{ (u,0): u\in \mathbb{R}^n \}$ and the corresponding action of $G$ on $C_0(X)$ is almost-free.

\smnoind
(c) $G=SO(2)$ acts on $X=\mathbb{R}^2$ by rotation. 
For any $t\in SO(2)\setminus \{e\}$, we have $X^t = \{0\}$ and again the corresponding action of $G$ on $C_0(X)$ is almost-free.

\smnoind
(d) Let $G$ be any amenable group and $X = G \cup \{x_0\}$ be the one point compatification of $G$. 
Define an action $\alpha$ of $G$ on $X$ by 
$$ \alpha_t(x) = 
\begin{cases}
tx  & \qquad {\rm if}{\ }x\in G\\
x_0  & \qquad {\rm if}{\ }x=x_0.
\end{cases}$$
It is not hard to see that this action is continuous. 
Moreover, for any $t\in G\setminus \{e\}$, we have $X^t = \{ x_0\}$ and $\alpha$ is almost-free.

\smnoind
(e) Let $f$ be the function from $\mathbb{R}$ to $\mathbb{C}$ defined by $f(s) = (1-e^{-\abs{s}})e^{{\rm i}s}$ (for $s\in \mathbb{R}$). 
Let $X = \overline{ \{f(s): s\in \mathbb{R} \}} = \{f(s): s\in \mathbb{R} \}\cup \mathbb{T}$ ($\mathbb{T}$ is the torus). 
Define an action $\alpha$ of $\mathbb{R}$ on $X$ by $$ \alpha_t(x) = 
\begin{cases}
f(t+r)  & \qquad {\rm if}{\ }x=f(r) {\ }{\rm for{\ }some}{\ }r \in \mathbb{R}\\
e^{{\rm i}t}x  & \qquad {\rm if}{\ }x\in \mathbb{T}.
\end{cases}$$
It is not hard to see that $\alpha$ is a continuous action. 
Moreover, we have 
$$X^t = 
\begin{cases}
\emptyset & \qquad {\rm if}{\ }t\neq 2n\pi{\ }{\rm for{\ }all}{\ }n\in \mathbb{Z}\\
\mathbb{T} & \qquad {\rm if}{\ }t=2n\pi{\ }{\rm for{\ }some}{\ } n\in\mathbb{Z}\setminus\{0\}. 
\end{cases}$$
As $\mathbb{T}$ is nowhere dense in $X$, we see that $\alpha$ is almost-free.

\smnoind
(f) If $(\alpha,u)$ is almost-free and $H$ is a closed subgroup of $G$, then $(\alpha\!\mid_H, u\!\mid_{H\times H})$ is also almost-free. 
\end{eg}

\medskip

The above are examples of almost-free actions that are not free. 
In the following, we will give some examples of topologically free actions that are not almost-free.
Moreover, the Connes' spectrum of the second example is the whole torus.

\begin{eg}
\label{non-eg}
(a) Suppose that $G$ is the discrete Heisenberg group 
$\left\{ \left(
\begin{array}{ccc}
  1 & i & j \\
  0 & 1 & k \\
  0 & 0 & 1 \\
\end{array}
\right): i,j,k \in \mathbb{Z}\right\}$
acting canonically on $X = \mathbb{R}^3$. 
Take any $t = \left(%
\begin{array}{ccc}
  1 & i & j \\
  0 & 1 & k \\
  0 & 0 & 1 \\
\end{array}%
\right) \in G\setminus \{e\}$. 
If $k\neq 0$, then $X^t$ is either $\mathbb{R}\times (0)\times (0)$ or $\mathbb{R}^2\times (0)$ (depending on whether $i=0$). 
If $k = 0$, then $X^t = \left\{ \left(%
\begin{array}{c}
  x \\
  y \\
  z \\
\end{array}%
\right)\in \mathbb{R}^3: iy= -jz\right\}.$
Hence, if $\beta$ is the corresponding action of $G$ on $A = C_0(X)$, then
$$\pf(A) = \left\{ \left(%
\begin{array}{c}
  x \\
  y \\
  z \\
\end{array}%
\right)\in \mathbb{R}^3: y = az \ {\rm for\ some}\ a\in \mathbb{R}\setminus \mathbb{Q}\right\}$$
which is dense in $X$ but $\I(\pf(A)) = \emptyset$.
This shows that $\beta$ is topologically free but not almost-free. 

\smnoind
(b) Let $X$ be the set of all functions from $\mathbb{Z}$ to $\{0,1\}$. 
We consider $X = \prod_{-\infty}^{\infty} \{0,1\}$ with the product topology. 
Let $\alpha$ be the continuous action of $\mathbb{Z}$ on $X$ defined by 
$$\alpha_n(f)(m) = f(n+m).$$ 
For any $n\in \mathbb{Z}\setminus \{0\}$, it is clear that $X^n$ is the set of all periodic functions of period $\abs{n}$. 
Hence $X^n$ does not contain any set in the canonical base $\cal B$ of the product topology on $X$ and the corresponding action $\beta$ of $\mathbb{Z}$ on $C(X)$ satisfies \cite[10.4(iv)]{OP3}. 
Now, \cite[10.4(v)]{OP3} shows that $\beta$ is topologically free. 
On the other hand, any set in $\cal B$ will have non-empty intersection with $\bigcup_{n\in \mathbb{Z}\setminus \{0\}} X^n$. 
This shows that $\bigcup_{n\in \mathbb{Z}\setminus \{0\}} X^n$ is dense in $X$ and $\beta$ is not almost-free. 
\end{eg}

\medskip

\begin{thm}
\label{twisted-case}
Let $(A,G,\beta,u)$ be a separable twisted $C^*$-dynamical system such that $G$ is amenable. 
If $(\beta, u)$ is almost-free, then any non-zero ideal $J_0$ of $A\times_{\beta,u} G$ contains an ideal of the form $I_0\times_{\beta,u} G$ where $I_0$ is a non-zero $\beta$-invariant ideal of $A$. 
\end{thm}
\noindent
{\bf Proof:}
We consider, first of all, the situation when $(\beta,u)$ is an ordinary action (i.e. $u$ is trivial and $\beta$ is continuous). 
Let $D = A\times_{\beta} G$. 
Let $W_0 = \{ J\in \pr(D): J_0\nsubseteq J\}\neq \emptyset$. 
Let $\Phi$ be the quasi-orbit map from $\pr(D)$ to $\pq(\beta)$.
Then $$\emptyset\neq V_0 = q_G^{-1}(\Phi(W_0))\subseteq \pr(A)$$ 
is an open set because $\Phi$ is open (\cite[4.8]{GL}). 
As $V_0$ is $G$-invariant, Proposition \ref{almost-free}(iv) implies that 
$$U_0 = \I(V_0\cap {\rm Prim}_\beta (A))$$ 
is an non-empty open $G$-invariant subset of $\pr(A)$.
Hence $$I_0 = \ker (\pr(A)\setminus U_0)$$ is a non-zero $\beta$-invariant ideal of $A$. 
It remains to show that $I_0\times_\beta G \subseteq J_0$, or equivalently, $i(x)j(h)\in J$ for any $x\in I_0$, $h\in C_c(G)$ and $J\in \pr(D)\setminus W_0$ (where $i$ and $j$ are the canonical maps from $A$ and $C^*(G)$ respectively to $M(D)$). 
Fix any $x\in I_0$, $h\in C_c(G)$ and $J\in \pr(D)\setminus W_0$. 
By \cite[3.2]{GR}, there exists $I\in \pr(A)$ such that $J$ lives on the quasi-orbit $q_G^{\ }(I)$ and there is a representation $\mu = \pi\times v$ (on a Hilbert space $\cal K$) of $A\times_\beta G\!_I$ (where $G\!_I$ is the stabilizer of $I$) such that 
$$\ker (\pi) = I \quad {\rm and} \quad J = \ker(\nu)$$ 
where $\nu = \ind_{G\!_I}^G(\mu)$. 
It suffices to show that $\nu(i(x)j(h)) = 0$. 
First of all, we claim that $I\notin U_0$. 
Suppose on the contrary that $I\in U_0 \subseteq V_0$.
Then, there exists $J'\in W_0$ such that $$\Phi(J') = q_G^{\ }(I)$$
(i.e. $J'$ lives on $q_G^{\ }(I)$). 
On the other hand, as $U_0\subseteq {\rm Prim}_\beta (A)$, we have $J'=J$ by Lemmas \ref{inv-open}(a) \& \ref{unique}. 
This contradicts the fact that $J\notin W_0$. 
Next, we recall that $\nu$ is a representation of $D$ on the Hilbert space $\cal H$ generated by $\{ F\otimes_\mu \xi: F\in C_c(G;A); \xi\in \mathcal{K} \}$ and is also generated by $\{f\otimes_\mu \xi: f\in C_c(G); \xi\in \mathcal{K} \}$. 
Moreover, $$\nu(i(x)j(h))(f\otimes_\mu \xi) = i(x)j(h*f)\otimes_{\mu} \xi$$ 
(where $h*f$ is the convolution product of $h$ and $f$ in $C_c(G)$). 
Note that for any $F\in C_c(G,I_0)$, we have 
$$\mu(F\!\mid_{G_I})\xi = \int_{G_I} \pi(F(s))v(s)\xi{\ }ds = 0$$ 
because $I_0\subseteq I = \ker(\pi)$ (as $I \in \pr(A)\setminus U_0$).
Thus, 
$$\| i(x)j(h*f)\otimes_\mu \xi\|^2 \, = \,\abs{\la\mu(j(h*f)^*i(x^*x)j(h*f)\!\mid_{G_I}) \xi, \xi \ra} \, = \, 0.$$ 
This shows that $\nu(i(x)j(h)) =0$ as required. 
Now, we turn to the twisted case.
We recall that by the Packer-Raeburn stabilisation trick in \cite[3.4]{PR}, there exists a continuous action $\gamma$ of $G$ on $B = A\otimes \mathcal{K}$ (where $\mathcal{K} = \mathcal{K}(L^2(G))$) such that $(\beta\otimes \id, u\otimes 1)$ is exterior equivalent to $\gamma$ in the sense of \cite[3.1]{PR}. 
Notice that $B$ is separable as both $A$ and $G$ are. 
Moreover, if $$\Theta: \pr(A) \rightarrow \pr(B)$$ 
is the canonical homeomorphism, then $G$ acts freely on an open subset $V\subset \pr(A)$ if and only if $G$ acts freely on $\Theta(V)$ (note that the induced action of $\beta\otimes \id$ on $\pr(B)$ and that of $\gamma$ are the same). 
This shows that $\gamma$ is also almost-free (by Proposition \ref{almost-free}).
Let $$\varphi: (A\times_{\beta, u} G) \otimes \mathcal{K} \rightarrow B\times_\gamma G$$ be the canonical isomorphism and $J' = \varphi(J_0\otimes \mathcal{K})$. 
By the untwisted case above, there exists a non-trivial $\gamma$-invariant ideal $I'$ of $B$ such that $I'\times_\gamma G \subseteq J'$. 
As $I'$ is also $\beta\otimes \id$-invariant, $I' = I_0 \otimes \mathcal{K}$ for an $\beta$-invariant ideal $I_0$. 
It is clear that $\varphi((I_0\times_{\beta, u} G)\otimes \mathcal{K}) = I'\times_\gamma G \subseteq \varphi(J_0\otimes \mathcal{K})$ and so $I_0\times_{\beta, u} G \subseteq J_0$. 

\medskip

\begin{rem}
The conclusion of Theorem \ref{twisted-case} means that for any representation $\nu$ of $A\times_{\beta,u} G$, if $\nu\!\mid_A$ is injective, then so is $\nu$. 
\end{rem}

\medskip

Moreover, we can also use this theorem to obtain the following result concerning the primeness (or primitivity) of a twisted crossed product.  

\medskip

\begin{cor}
\label{prime}
Under the assumption of Theorem \ref{twisted-case}, the following statements are equivalent. 

\smnoind
(1) $A$ is $G$-prime. 

\smnoind
(2) $A\times_{\beta,u} G$ is prime.

\smnoind
(3) There exists a dense $G$-orbit in $\pr (A)$. 

\smnoind
(4) $A\times_{\beta,u} G$ is primitive. 
\end{cor}
\noindent
{\bf Proof:}
Suppose that there exist two non-trivial ideals $J_1$ and $J_2$ of $A\times_{\beta, u} G$ such that $J_1\cdot J_2 = (0)$.
By Theorem \ref{twisted-case}, there exist non-zero $G$-invariant ideals $I_1$ and $I_2$ of $A$ such that 
$$I_i\times_{\beta, u} G \subseteq J_i$$ 
($i=1,2$). 
As $I_1\cdot I_2 = (0)$, we see that $A$ is not $G$-prime.
This shows that (1) implies (2). 
The converse of this implication is trivial. 
To show that (3) is equivalent to (4), we consider first of all, the case when $(\beta,u)$ is an ordinary action.  
Suppose that $A\times_{\beta} G$ is primitive and $\mathcal{Q} = \Phi((0))$ (where $\Phi$ is the quasi-orbit map as in Theorem \ref{twisted-case}).
Then $\cal Q$ is dense in $\pr(A)$ (as $\ker (\mathcal{Q}) = (0)$). 
Hence any $G$-orbit contained in $\mathcal{Q}$ is dense (as $\mathcal{Q}$, being a quasi-orbit, lies in the closure of any $G$-orbit that it contains). 
Conversely, suppose that $\mathcal{O}$ is a dense $G$-orbit and $\cal Q$ is the quasi-orbit containing $\mathcal{O}$. 
Then $\ker (\mathcal{Q}) = (0)$ (because $\cal Q$ is dense). 
By the surjectivity of $\Phi$ (\cite[4.8]{GL}), there exists $J\in \pr(A\times_{\beta} G)$ that lives on $\cal Q$. 
If $J\neq (0)$, then by Theorem \ref{twisted-case}, there is a non-trivial $G$-invariant ideal $I$ of $A$ such that $I\times_\beta G \subseteq J$ and thus, 
$$(0) = \ker (\mathcal{Q}) = \Res^G_{(e)}(J) \supseteq I$$ 
which implies a contradiction. 
Now, we can get the general case by using the Packer-Raeburn stabilisation trick as in Theorem \ref{twisted-case} (note that $A$ satisfies condition (3) if and only if $A\otimes \mathcal{K}$ does and $A\times_{\beta,u} G$ is primitive if and only if $(A\times_{\beta,u} G)\otimes \mathcal{K}$ is). 
Finally, as $A\times_{\beta,u} G$ is separable, it is primitive if and only if it is prime. 

\medskip

\begin{rem}
(a) It is known that in the separable case, a twisted covariant system in the sense of Green (see \cite[p.196]{Gre}) will induce a twisted $C^*$-dynamical system as above (see e.g. \cite[5.1]{PR}). 
Moreover, it is easy to see that one has the corresponding results for Theorem \ref{twisted-case} and Corollary \ref{prime} for twisted covariant systems (in this case, almost-freeness means that $\bigcup_{t\in G\setminus N} \pr(A)^t$ is nowhere dense in $\pr(A)$). 
Using this, one can show that if $G$ is a separable group with an abelian closed normal subgroup $N$ such that $\bigcup_{t\in G\setminus N} \hat N^t$ is nowhere dense in the dual group $\hat N$, then $G$ is a $C^*$-unique group (in the sense of \cite{Boi-unique}; see \cite{LN} for the details). 

\smnoind
(b) Note that if the action of $G$ on $A$ is quasi-regular (see \cite[p.223]{Gre}) and every quasi-orbit is regular (see \cite[p.223]{Gre}), then by using \cite[propostion 20]{Gre} instead of \cite[3.2]{GR}, one can remove the separabilities for both $G$ and $A$ in Lemma \ref{unique} as well as in Theorem \ref{twisted-case} and obtain certain untwisted versions of Theorem \ref{twisted-case} and Corollary \ref{prime} (in this case, we only have (1) being equivalent to (2) and (3) being equivalent to (4)).

\smnoind
(c) Note that in the case of an ordinary action, we can identify the $G$-primitive spectrum of $A$ (as defined in \cite{Mas}) with the quasi-orbit space since the quasi-orbit map is surjective and two quasi-orbits are the same if and only if they have the same kernel (note also that the topology given in \cite{Mas} is the same as the quotient topology on the quasi-orbit space). 
In this case, condition (3) of Corollary \ref{prime} holds if and only if $(0)$ is an element in the $G$-primitive spectrum of $A$ (by the argument of Corollary \ref{prime}). 
Thus, we say that $A$ is \emph{$G$-primitive} if condition (3) of Corollary \ref{prime} holds. 

\smnoind
(d) In the case of discrete groups (not necessarily countable), we can use the twisted version of \cite[theorem 1]{AS} (see \cite[A.1]{LN}) to give the following stronger version of Corollary \ref{prime}:

\smnoind
Suppose that $\Gamma$ is a discrete amenable group and $(\alpha, \tau)$ is a twisted action of $\Gamma$ on a $C^*$-algebra $B$ over a normal subgroup $\Delta$ (in the sense of \cite[p.196]{Gre}) such that $\bigcap_{i=1}^n \{[\pi]\in \hat B: t_i\cdot [\pi] \neq [\pi] \}$ is dense in $\hat B$ for any $t_1,...,t_n\in \Gamma\setminus \Delta$. 

\smnoind
(i) $B$ is $\Gamma$-prime if and only if $B\times_{\alpha,\tau} (\Gamma,\Delta)$ is prime. 

\smnoind
(ii) $B$ is $\Gamma$-primitive if and only if $B\times_{\alpha,\tau} (\Gamma,\Delta)$ is primitive. 

\smnoind
A similar result for twisted actions in the sense of Busby and Smith is also true (but separability will be needed because of the assumption in \cite{PR}).
\end{rem}

\medskip

We can prove a similar type of result for the simplicity of crossed products. 
In fact, we even have a better result (in the sense that $\I(\pf(A))$ is only required to be non-empty instead of dense) and do not need Theorem \ref{twisted-case}. 

\medskip

\begin{prop}
Let $(A,G,\beta,u)$ be a separable twisted $C^*$-dynamical system such that $G$ is amenable and $\I(\pf(A))$ is non-empty. 
Then $A$ is $G$-simple if and only if $A\times_{\beta,u} G$ is simple. 
\end{prop}
\noindent
{\bf Proof:}
We consider the case of ordinary action first. 
Let $I = \ker(\bigcup_{t\in G\setminus \{e\}} \pr(A)^t)$. 
Then $\I(\pf(A)) = \{ J\in \pr(A): I\nsubseteq J\}$. 
By the $G$-simplicity, $I=A$ and therefore, $\pf(A)=\pr(A)$. 
Now, we can apply \cite[3.3]{GR} to conclude that $A\times_{\beta} G$ is simple.
The converse is clear and the twisted case follows again from the stabilisation trick as in Theorem \ref{twisted-case}. 

\bigskip

We end this paper with the following simple results concerning the simplicity of a crossed product as well as a discussion on the strong Connes' spectrum (for ordinary actions). 
Let us first recall in the following the canonical coaction for twisted crossed products (in Green's sense). 
From now on, $N$ is a closed normal subgroup of $G$ (not necessary separable). 

\medskip

\begin{rem}
(a) Suppose that $(\gamma, \tau)$ is a twisted action of $G$ on a $C^*$-algebra $B$ over $N$ (in the sense of \cite[p.196]{Gre}). 
Let $Q_B$ be the canonical quotient map from $B\times_\gamma G$ to $B\times_{\gamma, \tau} (G,N)$ and $I = \ker (Q_B)$.
We identify $B$ and $C^*(G)$ as subalgebras of $M(B\times_\gamma G)$. 
For any $n\in N$, let $\lambda_n$ be the canonical image of $n$ in $M(C^*(G))$. 
Then clearly, $Q_B(\lambda_n) = Q_B(\tau(n))$. 
Suppose that $(\pi, v)$ is a covariant representation for $(B, G, \gamma)$ such that $$\ker ((Q_B\otimes Q_N)\circ \hat \gamma) \subseteq \ker (\pi\times v)$$ 
(where $Q_N$ is the canonical quotient map from $C^*(G)$ to $C^*(G/N)$ and $\hat \gamma$ is the dual coaction). 
Then $\ker (\overline{(Q_B\otimes Q_N)\circ \hat \gamma}) \subseteq \ker (\overline{\pi\times v})$ (where $\overline{(Q_B\otimes Q_N)\circ \hat \gamma}$ and $\overline{\pi\times v}$ are the corresponding extensions to $M(B\times_\gamma G)$). 
Since $\overline{(Q_B\otimes Q_N)\hat \gamma}(\lambda_n - \tau(n)) = (Q_B\otimes \id)(\lambda_n\otimes 1 - \tau(n) \otimes 1) = 0$, we have 
$$v(n) - \pi(\tau(n)) = \overline{\pi\times v}(\lambda_n - \tau(n)) = 0.$$ 
This shows that $I\subseteq \ker ((Q_B\otimes Q_N)\circ \hat \gamma)$ (as it is the smallest ideal for which the above condition holds) and we have a map $\hat \gamma_\tau: B\times_{\gamma, \tau} (G,N) \rightarrow M(B\times_{\gamma, \tau} (G,N)\otimes_{\rm max} C^*(G/N))$ such that $\hat \gamma_\tau \circ Q_B = (Q_B\otimes Q_N)\circ \hat \gamma$. 
This is called the \emph{dual coaction} of $(\gamma, \tau)$. 

\smnoind
(b) For any (full) coaction $\delta: A\rightarrow M(A\otimes_{\rm max} C^*(G))$ of $G$ on $A$, one can define the restriction $\delta\!\mid_{G/N} = (\id\otimes Q_N)\circ \delta$ which is a coaction of $G/N$ (see e.g. \cite{Rae}). 

\end{rem}

\medskip

We recall that for any ordinary action $\alpha$ of $G$ on $A$, there is a twisted action $(\alpha^N, \tau)$ of $G$ on $A\times_{\alpha} N$ over $N$ such that there exists a $*$-isomorphism $\Psi: A\times_\alpha G \rightarrow (A\times_{\alpha} N)\times_{\alpha^N, \tau} (G,N)$ (see \cite[proposition 1]{Gre}). 

\medskip

\begin{lem}
An ideal $J$ of $A\times_{\alpha} G$ is $\hat\alpha\!\mid_{G/N}$-invariant (where $\hat \alpha$ is the dual coaction of $\alpha$) if and only if there is an $\alpha^N$-invariant ideal $I$ of $A\times_{\alpha} N$ such that $\Psi(J) = I\times_{\alpha^N, \tau} (G,N)$.
\end{lem}
\noindent
{\bf Proof:}
We need to show that $\widehat{\alpha^N}\circ \Psi = (\Psi\otimes Q_N)\circ \hat \alpha$. 
In fact, it is clear that the equality holds on the canonical images of both $A$ and $G$ in $M(A\times_\alpha G)$ and the lemma follows from a simple continuity argument. 

\bigskip

The following proposition is now a direct consequence of this lemma. 

\medskip

\begin{prop}
Let $(A,G,\alpha)$ be a $C^*$-dynamical system and $N$ is any closed normal subgroup of $G$. 
Then $A\times_\alpha G$ is simple if and only if $A\times_\alpha N$ is $\alpha^N$-simple and for any $J\in \pr(A\times_\alpha G)$, $J$ is $\hat \alpha\!\mid_{G/N}$-invariant. 
\end{prop}

\medskip

\begin{rem}
Suppose that $G$ is abelian and $\Gamma$ is the dual group of $G$. 

\smnoind
(a) One equivalent form of the strong Connes' spectrum is that $\ti\Gamma (\alpha) = \{ \gamma\in \Gamma: \hat\alpha_\gamma(J)\subseteq J$ for any ideal $J$ of $A\times_{\alpha} G\}$ (see e.g. \cite[p.367]{Kis}).  
Therefore, we can reformulate the above proposition as follows (which is a strong Connes' spectrum version of \cite[3.1]{Ole}): 

\begin{equation}
\label{quote}
A\times_\alpha G {\rm\ is\ simple\ if\ and\ only\ if\ } A\times_\alpha N {\rm\ is\ } \alpha^N{\rm -simple\ and\ } N^\bot \subseteq \ti\Gamma(\alpha) 
\end{equation}
(note that $N^\bot$ is the image of $\widehat{G/N}$ in $\hat G$).

\smnoind
(b) Suppose that $H$ and $K$ are two subgroups of the semigroup $\ti\Gamma(\alpha)$. 
It is easy to see that the sub-semigroup generated by $H$ and $K$ is again a subgroup of $\ti\Gamma(\alpha)$. 
We denote by $\ti\Gamma_{\rm sym}(\alpha)$ the union of all subgroups in $\ti\Gamma(\alpha)$. 
Then $\ti\Gamma_{\rm sym}(\alpha)$ is a closed subgroup of $\Gamma$ (as $\ti\Gamma(\alpha)$ is closed) and we can replace the last part of the equivalent statement (\ref{quote}) by ``$N^\bot \subseteq \ti\Gamma_{\rm sym}(\alpha)$''. 

\smnoind
(c) In the case when $G$ is non-abelian, the ``non-commutative space'' $\ti\Gamma(\alpha)$ could be difficult to describe. 
However, we can identified a ``non-commutative closed subgroup of $\ti\Gamma(\alpha)$'' with its dual group (which is a quotient group of $G$). 
More precisely, $G/N$ is the dual group of such a ``subgroup'' if and only if any $J\in \pr(A\times_\alpha G)$ is $\hat \alpha\!\mid_{G/N}$-invariant. 
In this way, we can vaguely describe the ``non-commutative space'' $\ti\Gamma_{\rm sym} (\alpha)$. 
Of course, it will be nice if one can replace the description in terms of $A$ and $\alpha$ alone without considering $\pr(A\times_\alpha G)$. 
However, even in the abelian case, the original description of $\ti\Gamma(\alpha)$ (see \cite[p.366]{Kis}) depends on the space of irreducible representations of $A\times_\alpha G$. 
\end{rem}

\bigskip

\smnoind
Chi-Wai Leung, Department of Mathematics, The Chinese University of Hong Kong, Shatin, Hong Kong. 

\smnoind
\emph{Email address:} cwleung@math.cuhk.edu.hk

\bigskip\noindent
Chi-Keung Ng, Department of Pure Mathematics, Queen's University Belfast, Belfast BT7 1NN, United Kingdom. 

\smnoind
\emph{Email address:} c.k.ng@qub.ac.uk

\end{document}